\documentclass[fleqn]{article}
\usepackage{amsmath, amssymb, amsthm, cite}
\usepackage[dvipdfm]{graphicx}

\usepackage[dvipdfm]{color}
\usepackage[dvipdfm,bookmarks=true,bookmarksnumbered=true,bookmarkstype=toc]{hyperref} 
\newtheorem{Thm}{Theorem}[section]
\newtheorem{Def}[Thm]{Definition}
\newtheorem{Lem}[Thm]{Lemma}
\newtheorem{Prop}[Thm]{Proposition}
\newtheorem{Cor}[Thm]{Corollary}
\newtheorem{Rem}[Thm]{Remark}
\newcommand{\myfig}[3]{%
\begin{figure}[htbp]%
\centering \includegraphics{#1}%
\caption{#2} \label{#3}%
\end{figure}}
\begin{document}

\title{Combinatorial proofs for basic properties of Ozsv\'ath-Szab\'o invariant}
\author{Shinya Ichida}
\maketitle

\begin{abstract}
Sarkar and Wang have given a combinatorial algorithm for computing Heegaard Floer homology and Plamenevskaya has improved their method to compute Ozsv\'ath-Szab\'o invariant.

In this paper, applying the combinatorial method to stabilizations of an open book, we prove basic properties of Ozsv\'ath-Szab\'o invariant.
\end{abstract}

\section{Introduction} \label{Intro}
Let $(Y, \xi)$ be a $3$-dimensional contact manifold. $(Y, \xi)$ is called \textit{overtwisted}, if $(Y, \xi)$ has an embedded disk $D$ such that $\partial D$ is tangent to $\xi$ and the framing of $\partial D$ given by $\xi$ coincides with the framing given by $D$. $(Y, \xi)$ is called \textit{tight}, if it is not overtwisted. Eliashberg has classified overtwisted contact structures by using homotopy theory in \cite{YE1989}. However, tight contact structures have not been classified completely. 

One of basic theorems in $3$-dimensional contact topology is Giroux's correspondence. It says that there exists a one-to-one correspondence between isotopy equivalence classes of contact structures and positive stably equivalence classes of open books.

Using Giroux's correspondence, Ozsv\'ath and Szab\'o have defined an invariant of isotopy classes of contact structures in \cite{OS2005}. The invariant is an element of Heegaard Floer homology of the ambient manifold. It is useful for classifying contact structures up to isotopy. For example, if the contact structure is overtwisted, then their invariant vanishes. In particular, it is also useful for classifying tight contact structures.

Sarkar and Wang have given a combinatorial algorithm for computing Heegaard Floer homology in \cite{SW2008} and Plamenevskaya has improved their method to compute Ozsv\'ath-Szab\'o invariant in \cite{OP2007}.

In this paper, after reviewing their results, applying the combinatorial method to stabilizations of an open book, we prove the following two basic properties of Ozsv\'ath-Szab\'o invariant.

Let $\xi$ be a contact structure on a $3$-manifold, and $c( \xi )$ be the Ozsv\'ath-Szab\'o invariant of $\xi$.
\begin{itemize}
 \item $c( \xi )$ vanishes, if $\xi$ is compatible with a negative stabilization of any open book.
 \item Positive stabilization preserves $c( \xi )$.
\end{itemize}
This paper also includes an easier computation of the contact invariant of the example in \cite{OP2007}.

This paper is the author's master's thesis. I would like to thank my supervisor, Masaaki Ue. I would also like to thank Ryokichi Tanaka. I could not write this paper without their adivice and encouragement.

\section{Giroux's correspondence} \label{Giroux-corresp}
In this section, we review Giroux's correspondence theorem.

An \textit{open book} is a pair $(S, h)$, where $S$ is the \textit{page}, which is an oriented surface with non-empty boundary, and $h$ is the \textit{monodromy}, which is an orientation preserving self-diffeomorphism of $S$ which coincides with the identity map on $\partial S$.  We call two open books $(S_{1}, h_{1})$ and $(S_{2}, h_{2})$ are equivalent, if $S_{1}$ equals $S_{2}$ and $h_{1}$ is isotopic to $g \circ h_{2} \circ g^{-1}$ relative $\partial S$, where $g$ is an orientation preserving self-diffeomorphism of $S$.

Let $(Y, \xi)$ be a contact $3$-manifold. $(Y, \xi)$ and an open book $(S, h)$ are called compatible, if $Y$ is diffeomorphic to $S \times [0, 1] / \sim$, where the equivalence relation $\sim$ is given by
\[ (x, 1) \sim (h(x), 0), \qquad x \in S \]
\[ (x, t) \sim (x, t'), \qquad x \in \partial S, \, t, t' \in [0, 1] \]
and the following conditions hold.
\begin{itemize}
  \item $\partial S$ is a transverse link, that is, its nonzero tangent vector is transverse to $\xi$,
  \item $d \alpha$ is a volume form on $S$, where $\alpha$ is a contact form of $\xi$,
  \item the orientation of $\partial S$ which satisfies $\alpha > 0$ coincides with the boundary orientation induced by the surface orientation defined by $d \alpha$.
\end{itemize}

For an open book $(S, h)$, its \textit{positive(negative) stabilization} is an open book $(S', h')$ such that
\begin{itemize}
  \item $S'$ is a surface obtained from $S$ by attaching a $1$-handle to $S$,
  \item $h'$ is a composition map $h' = t_{\gamma} \circ h $, where
  \item $\gamma$ is a simple closed curve in $S'$ which goes over the $1$-handle once, and $t_{\gamma}$ is the positive(negative) Dehn twist around it. 
\end{itemize}
This operation does not change the ambient manifold $Y$. Note that a positive(negative) stabilization does depend on $\gamma$.

We call two open books $(S_{1}, h_{1})$ and $(S_{2}, h_{2})$ are positive stably equivalent, if a finitely many positive stabilizations of $(S_{1}, h_{1})$ is equivalent to a finitely many positive stabilizations of $(S_{2}, h_{2})$. 

Giroux has shown the following theorem in \cite{EG2002}.
\begin{Thm} \label{Giroux}
Let $(Y, \xi)$ be a $3$-dimensional contact manifold. Then there exists a unique compatible open book $(S, h)$ up to positive stably equivalence relation. Conversely, let $(S, h)$ be an open book. Then there exists a unique compatible contact manifold $(Y, \xi)$ up to isotopy.
In addition, this correspondence is one-to-one.
\end{Thm}
\begin{Rem}
A positive stabilization does not change the compatible contact structure. However, this is not the case for a negative stabilization. It is known that a contact structure compatible with a negative stabilization of any open book is overtwisted, even if the original contact structure is tight.
\end{Rem}

\section {Hat Heegaard Floer homology}
In this section, we review the definition of hat Heegaard Floer homology, following \cite{SW2008}. 

Originally, hat Heegaard Floer homology $\widehat {HF}$ of a $3$-manifold is defined by using the moduli space of holomorphic representations of Whitney disks in the symmetric product of $\Sigma$, which is a Heegaard surface for a Heegaard decomposition of the $3$-manifold (see \cite{OS2004} and the review of the construction in \cite{OS2006}). In this paper, however, we use cylindrical reformulation of $\widehat {HF}$ by Lipshitz in \cite{RL2006}.

Let $Y$ be a closed oriented $3$-manifold and $Y = H_{1} \cup _{\Sigma} H_{2}$ be its Heegaard decomposition. $\Sigma = \partial H_{1} = \partial H_{2}$ is oriented as the boundary of $H_{1}$.

Let $g$ be the genus of $\Sigma$. Then, $Y$ has a handle decomposition having one $0$-handle, $g$ $1$-handles, $g$ $2$-handles, and one $3$-handle. We call the cocores of the $1$-handles $\alpha _{1} , \alpha _{2} , \ldots , \alpha _{g}$ and the attaching circles of the $2$-handles $\beta _{1} , \beta _{2} , \ldots , \beta _{g}$.  Then, we get a tuple $(\Sigma , \boldsymbol \alpha = ( \alpha _{1} , \alpha _{2} , \ldots , \alpha _{g} ) , \boldsymbol \beta  = ( \beta _{1} , \beta _{2} , \ldots , \beta _{g} ) )$ called a \textit{Heegaard diagram}. 

We fix a point $z_{0} \in \Sigma \setminus \left ( \bigcup \alpha _{i} \right ) \setminus \left ( \bigcup \beta _{j} \right )$ called a \textit{basepoint}. The tuple $(\Sigma , \boldsymbol { \alpha } , \boldsymbol { \beta } , z_{0} )$ is called a \textit{pointed Heegaard diagram}.

Given a pointed Heegaard diagram $(\Sigma , \boldsymbol { \alpha } , \boldsymbol { \beta } , z_{0} )$, the chain complex $\widehat{CF}$ is a free abelian group with $\mathbb{Z}/2\mathbb{Z}$-coefficients whose generators are formal sums of $g$ distinct points in $\Sigma$, $\mathbf{x} = x_{1} + x_{2} + \cdots + x_{g}$, such that each $\alpha$-curve contains some $x_{i}$ and each $\beta$-curve contains some $x_{j}$. That is, there exist $\sigma, \tau \in S_{g}$ such that $\alpha _{\sigma (i)} \cap \beta _{\tau (i)}$ includes $x_{i}$ for each $i = 1, 2, \ldots, g$. 
A connected component of $\Sigma \setminus (\boldsymbol{\alpha} \cup \boldsymbol{\beta})$ is called a \textit{region}. The region containing $z_{0}$ is called the \textit{pointed region}. A formal sum of regions with integer coefficients is called a $2$\textit{-chain}.
For a $2$-chain $\phi = a_{1} D_{1} + a_{2} D_{2} + \cdots + a_{r} D_{r}$, $\partial ( { \partial ( \phi ) } | _{\boldsymbol \alpha})$ is defined by
\[ \partial ( { \partial ( \phi ) } | _{\boldsymbol \alpha}) = a_{1} \partial ( { \partial ( D_{1} ) } | _{\boldsymbol \alpha}) + a_{2} \partial ( { \partial ( D_{2} ) } | _{\boldsymbol \alpha}) + \cdots + a_{r} \partial ( { \partial ( D_{r} ) } | _{\boldsymbol \alpha}), \]
where $\partial ( { \partial ( D ) } | _{\boldsymbol \alpha})$ is a formal sum of the endpoints of $\alpha$-edges of $D$ with signature which is induced by the orientation of $\Sigma$.
Given two generators $\mathbf{x}$, $\mathbf{y}$ of $\widehat{CF}$, we define $\pi _{2} ( \mathbf{x}, \mathbf{y} )$ to be the collection of all $2$-chains $\phi$ such that $\partial ( { \partial ( \phi ) } | _{\boldsymbol \alpha}) = \mathbf{y} - \mathbf{x}$. Such $2$-chains are called \textit{domains} connecting $\mathbf{x}$ to $\mathbf{y}$.
Given $p \in \Sigma \setminus ( \boldsymbol{\alpha} \cup \boldsymbol {\beta} )$ and $2$-chain $\phi$, let $n_{p} (\phi)$ be the coefficient of the region containing $p$ in $\phi$. A domain $\phi$ is \textit{positive}, if $n_{p}(\phi) \geq 0$ for every point $p$ in $\Sigma \setminus ( \boldsymbol{\alpha} \cup \boldsymbol {\beta} )$. We define 
$\pi _{2} ^{0} ( \mathbf{x}, \mathbf{y} ) = \{ \phi \in \pi _{2} ( \mathbf{x}, \mathbf{y} ); \; n_{ z_{0} } ( \phi ) = 0 \}$. A pointed Heegaard diagram is \textit{admissible}, if any positive domain $\phi \in \pi _{2} ^{0} ( \mathbf{x}, \mathbf{x} )$ is trivial, for every generator $\mathbf{x}$. It is known that there exists an admissible pointed Heegaard diagram for any $Y$. Hence, we assume that $(\Sigma , \boldsymbol { \alpha } , \boldsymbol { \beta } , z_{0} )$ is admissible.

For a domain $\phi$ and generator $\mathbf{x} = x_{1} + x_{2} + \cdots + x_{g}$, $\mu _{x_{i}} ( \phi )$ is defined to be the average of the coefficients of the four regions around $x_{i}$ in $\phi$. The \textit{point measure} $\mu _{\mathbf{x}} ( \phi )$ is defined as $\sum \mu _{x_{i}} ( \phi )$.

Any region has even vertices, because each of its vertices is intersection point between some $\alpha$-curve and some $\beta$-curve. Hence the boundary of the region consists of $\alpha$-edges and $\beta$-edges, which are placed alternately. If a region $D$ is a $2n$-gon, then we define $e(D) = 1 - \frac{n}{2}$, and for a domain $\phi = \sum a_{i} D_{i}$, its \textit{Euler measure} $e( \phi )$ is defined by $e(\phi) = \sum a_{i} e(D_{i})$.

Let $\mathcal{M} (\phi)$ be the moduli space of holomorphic representatives of $\phi \in \pi_{2} ( \mathbf{x} , \mathbf{y} )$, and $\mu (\phi)$ be the \textit{Maslov index} of $\phi$, which is an expected dimension of $\mathcal{M} (\phi)$. 
\begin{Prop}[Lipshitz \cite{RL2006}]
For a domain $\phi \in \pi_{2} ( \mathbf{x} , \mathbf{y} )$, the Maslov index of $\phi$ is given by 
\[ \mu ( \phi ) = e( \phi ) + \mu _{ \mathbf{ x } } ( \phi ) + \mu _{ \mathbf{ y } } ( \phi ). \]
\end{Prop}
If $\phi$ is non-trivial, $\mathcal{M} (\phi)$ admits a free $\mathbb{R}$ action. 
If $\mu ( \phi ) = 1$, $ \mathcal{ M } ( \phi ) / \mathbb{ R } $ is a zero-dimensional manifold, so we define its \textit{count function} $c( \phi )$ to be the number of points in $ \mathcal{ M } ( \phi ) / \mathbb{ R } $, counted modulo $2$.

After all, the boundary map of the chain complex $\widehat{CF}$ is given by
\[ \hat{\partial} \mathbf{x} = \sum _{ \mathbf{y} }
\sum _{ \substack{ \phi \in \pi _{2} ^{0} ( \mathbf{x} , \mathbf{y} ) \\
                   \mu (\phi) = 1 } }
  c( \phi ) \mathbf{y}, \]
for every generator $\mathbf{x}$. If the pointed Heegaard diagram is admissible, the right hand side is a finite sum. This map satisfies $\hat{\partial} \circ \hat{\partial} = 0$.
\begin{Thm}[Lipshitz \cite{RL2006}]
For a $3$-manifold $Y$, the homology of the chain complex $( \widehat{CF} , \hat{\partial} )$ is isomorphic to Ozsv\'ath-Szab\'o Heegaard Floer homology $\widehat{HF} ( Y )$, which is an invariant for the $3$-manifold $Y$.
\end{Thm}
Note that the only non-combinatorial part of this formulation is the count function $c( \phi )$.

\section{A combinatorial description of Heegaard Floer homology and Ozsv\'ath-Szab\'o contact invariant}
In this section, we review a combinatorial description of Heegaard Floer homology and Ozsv\'ath-Szab\'o contact invariant.

At first, we review a combinatorial description of Heegaard Floer homology by Sakar and Wang, in \cite{SW2008}.
\begin{Def}
A pointed Heegaard diagram $( \Sigma , \boldsymbol{ \alpha } , \boldsymbol { \beta } , z_{0} )$ is nice, if any region not containing $z_{0}$ is either a bigon or a square.
\end{Def}
It is known that a nice pointed Heegaard diagram is admissible(\cite{LMW2008}).
\begin{Def}
Let $( \Sigma , \boldsymbol{ \alpha } , \boldsymbol { \beta } , z_{0} )$ be a nice pointed Heegaard diagram.

A domain $\phi \in \pi _{2} ^{0} ( \mathbf{x} , \mathbf{y} )$ is an empty embedded $2n$-gon, if it satisfies the following conditions.
  \begin{itemize}
    \item The coefficient of each region in $\phi$ is $0$ or $1$.
    \item The closure of $\phi$ is topologically an embedded disk with $2n$ vertices on its boundary for some $n \in \mathbb{N}$, such that $\mu_{v} (\phi) = \frac{1}{4}$ for each vertex $v$. 
    \item The closure of $\phi$ does not contain any $x_{i}$ or $y_{i}$ in its interior, where $\mathbf{x} = \sum x_{i}$ and $\mathbf{y} = \sum y_{i}$.
  \end{itemize}
\end{Def}
\begin{Thm}[Sarkar, Wang \cite{SW2008}]
Let $( \Sigma , \boldsymbol{ \alpha } , \boldsymbol { \beta } , z_{0} )$ be a nice pointed Heegaard diagram.
  \begin{enumerate}
    \item If $\phi \in \pi _{2} ^{0} ( \mathbf{x} , \mathbf{y} )$ is  an empty embedded bigon or square, then $\mu (\phi) = 1$.
    \item Let $\phi$ be a domain such that $\mu ( \phi ) = 1$, then $c( \phi ) = 1$ if and only if  $\phi$ is either an empty embedded bigon or square. 
  \end{enumerate}
\end{Thm}
Hence, the boundary map $\hat{\partial} \colon \widehat{ CF } \rightarrow \widehat { CF }$ can be computed as
\[ \hat{\partial} \mathbf{x} = \sum _{ \mathbf{y} }
\sum _{ \substack{
            \phi \in \pi _{2} ^{0} ( \mathbf{x} , \mathbf{y} ) \\
            \phi \text { is an empty embedded} \\
            \text{ bigon or square } 
        } 
      }
    \mathbf{y} ,\]
for a nice diagram.
\begin{Thm}[Sarkar, Wang \cite{SW2008}]
For any admissible pointed Heegaard diagram $( \Sigma , \boldsymbol{ \alpha } , \boldsymbol { \beta } , z_{0} )$, we can get a nice diagram by finite isotopies of $\boldsymbol{ \alpha }$ and handle slides among $\boldsymbol{ \alpha }$. We can also get a nice diagram by finite isotopies of $\boldsymbol{ \beta }$ and handle slides among $\boldsymbol{ \beta }$.
\end{Thm}
The above two theorems show that hat Heegaard Floer homology of any $3$-manifold can be computed combinatorially.

Next, we describe a construction of Ozsv\'ath-Szab\'o invariant by Honda, Kazez, and Mati\'c, in \cite{HKM2007}. 

Let $( Y , \xi )$ be a contact $3$-manifold and $(S , h)$ be a compatible open book. So Y is diffeomorphic to $S \times [0, 1] / \sim$, where the equivalence relation $\sim$ is defined as in section \ref{Giroux-corresp}.
Then, we get a Heegaard decomposition
\[ Y = H_{1} \cup H_{2} \]
with $H_{1} = { S \times [0, 1/2] } / \sim$, $H_{2} = { S \times [1/2, 1]} / \sim$.

There is a construction of the Heegaard diagram corresponding to the decomposition. $\Sigma$ denotes $ \partial H_{1} = S_{1/2} \cup - S_{0} $, where $S_{0} = S \times \{ 0 \}$, $S_{1/2} = S \times \{ 1/2 \}$, and $\alpha$- and $\beta$- curves are defined as follows.

Let $\{ a_{1}, a_{2}, \ldots , a_{n} \}$, called a \textit{basis} for $S$, be a collection of disjoint properly embedded arcs in $S$ such that $S \setminus \left ( \bigcup a_{i} \right )$ is a single polygon. Then, $b_{i}$ is determined as an arc which is isotopic to $a_{i}$ by a small isotopy so that the following conditions hold.
\begin{itemize}
  \item The endpoints of $a_{i}$ is isotopied along $\partial S$, in the direction induced by the orientation of $S$.
  \item $a_{i}$ and $b_{i}$ intersect transversely at only one point in the interior of $S$.
  \item If we orient $a_{i}$, and $b_{i}$ is given the orientation induced from the isotopy, then the sign of the intersection $a_{i} \cap b_{i}$ is positive.
\end{itemize}
For examples of such $\{ a_{1}, a_{2}, \ldots , a_{n} \}$ and $\{ b_{1}, b_{2}, \ldots , b_{n} \}$, see  figures in section \ref{exam}.
Then, the curves $ \alpha _{i} = a_{i} \times \{ 1/2 \} \cup a_{i} \times \{ 0 \} $ are cocores of $1$-handles and
$ \beta _{i} = b_{i} \times \{ 1/2 \} \cup h( b_{i} ) \times \{ 0 \} $
 are attaching circles of $2$-handles. Hence, we can consider $\boldsymbol{ \alpha } = ( \alpha _{1}, \alpha _{2}, \ldots, \alpha _{n} )$ and $\boldsymbol{ \beta } = ( \beta _{1}, \beta _{2}, \ldots, \beta _{n} )$ as $\alpha$- and $\beta$- curves in the Heegaard diagram.

In addition, the basepoint $z_{0}$ is placed on $S_{1/2}$ and is not in any thin strips between $a_{i} \times \{ 1/2 \}$ and $b_{i} \times \{ 1/2 \}$.

After all, we get the pointed Heegaard diagram $( \Sigma, \boldsymbol{ \alpha }, \boldsymbol{ \beta }, z_{0} )$ given by the open book. This diagram depends on a monodromy $h$ and a basis for $S$. It is known that this diagram is admissible(\cite{HKM2007}).

Ozsv\'ath-Szab\'o invariant is an element of $\widehat{HF} (-Y)$, so we consider $-Y$, whose diagram corresponding to $-Y = H_{2} \cup _{\Sigma} H_{1}$ is given by $( \Sigma, \boldsymbol{ \beta }, \boldsymbol{ \alpha }, z_{0} )$, interchanging the roles of the $\alpha$- and $\beta$- curves.

Let $c_{i}$ be the intersection point between $a_{i} \times \{ 1/2 \}$ and $b_{i} \times \{ 1/2 \}$, for $i = 1, 2, \ldots, n$. A formal sum $\mathbf{c} = c_{1} + c_{2} + \cdots + c_{n}$ is a generator of $\widehat{ CF } ( \Sigma, \boldsymbol{ \beta }, \boldsymbol{ \alpha }, z_{0} )$.
\begin{Thm}[Honda, Kazez, Mati\'c \cite{HKM2007}] \label{HKM-th}
$\mathbf{c}$ is a cycle, which represents Ozsv\'ath-Szab\'o invariant $c( \xi )$ in $\widehat{ HF } (-Y)$.
\end{Thm}
Now, as an application of Sakar-Wang algorithm, we have the following theorem of Plamenevskaya.
\begin{Thm}[Plamenevskaya \cite{OP2007}] \label{Pla-th}
For any open book $(S, h)$ and any basis $\{ a_{1}, a_{2}, \ldots , a_{n} \}$ for $S$, there exists an orientation preserving diffeomorphism $\phi \colon S \rightarrow S$ which is isotopic to $\mathrm{id} _{S}$ relative $\partial S$ such that the open book $(S, \phi \circ h)$, which is equivalent to $(S, h)$, and the basis $\{ a_{1}, a_{2}, \ldots , a_{n} \}$ give a nice pointed Heegaard diagram.
\end{Thm}
We can get a combinatorial method of computing Ozsv\'ath-Szab\'o invariant by combining theorem \ref{HKM-th} and \ref{Pla-th}.
\begin{proof}[Sketch of proof of theorem \ref{Pla-th}]
Let $( \Sigma, \boldsymbol{ \alpha }, \boldsymbol{ \beta }, z_{0} )$ be the diagram given by $(S, h)$ and the basis $\{ a_{1}, a_{2}, \ldots , a_{n} \}$. We can get a nice diagram by moving $\beta$-curves of the original diagram by Sakar-Wang algorithm. If the moves include only isotopies in $S_{0} \subset \Sigma$, we prove the theorem.

There are two steps in Sakar-Wang algorithm. Firstly, we eliminate non-disk regions, then, make all but the pointed region bigons or squares.

\noindent
\textbf{Step.1 Eliminating the non-disk regions except the pointed region.}

Since $\Sigma \setminus \boldsymbol {\alpha}$ is a punctured sphere, each non-disk region is a punctured sphere with more than one boundary components. In addition, each boundary component dose contain both $\alpha$-edges and $\beta$-edges, because every $\alpha _{i}$ intersects $\beta _{i}$ in our diagram.

So, if $D$ is a non-disk region, then we can take disjoint, properly embedded arcs $l_{1}, l_{2}, \ldots, l_{k}$ in $D$, such that they are contained in $S_{0}$ and each $l_{i}$ connects an $\alpha$-curve to a $\beta$-curve, so that $D \setminus l_{1} \setminus \cdots \setminus l_{k}$ is a polygon.

Then, we make isotopies, called finger moves, by pushing $\beta$-curves along $l_{i}$ as in Figure \ref{zu-non-disk}. They are isotopies in $S_{0}$ which kill the non-disk regions.
\myfig{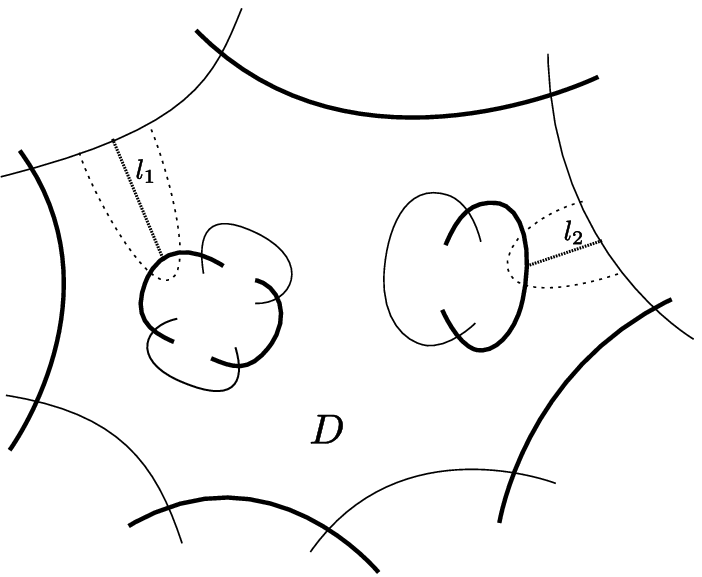}{Eliminating a non-disk region $D$. Here, thick solid arcs denote $\alpha$-curves and thin solid arcs denote $\beta$-curves. The same notation will be used in Figure \ref{zu-making-good}.}{zu-non-disk}

\noindent
\textbf{Step.2 Making all but the pointed region bigons or squares.}

Let $D_{0}$ denote the pointed region. We define the \textit{distance} $d(D)$ of a region $D$ to be the smallest number of intersection points between $\beta$-curves and an arc connecting $z_{0}$ to an interior point of $D$ in $\Sigma \setminus \boldsymbol{ \alpha }$.

We say a region $D$ is \textit{bad}, if $D$ is not pointed, and is not a bigon nor a square region. We take a bad region $D$ with the largest distance $d(D)$, and replace it by square regions. One of the regions having a common $\beta$-edge with $D$ has distance $d(D) - 1$ by definition. Let $D_{\ast}$ be such a region and $b_{\ast}$ be its common $\beta$-edge with $D$.

We make a finger move of $b_{\ast}$ into an adjacent region which has a common $\alpha$-edge with $D$ as in Figure \ref{zu-making-good} so that $D$ is replaced by two regions of fewer edges.
\myfig{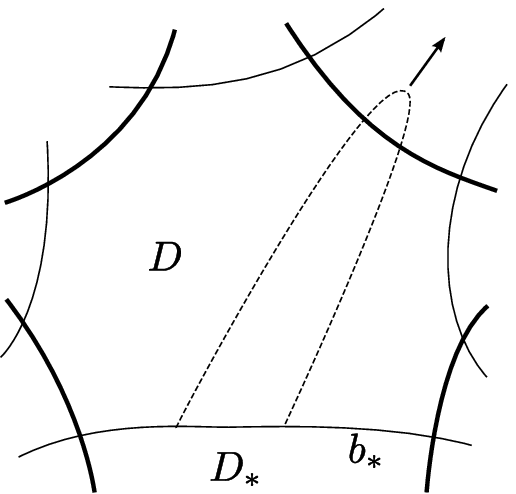}{Killing a bad region.}{zu-making-good}
The finger move is continued until the finger reaches either a bigon, another bad region, or a region of smaller distance. Plamenevskaya shows that, by further finger moves, the finger reaches either a bigon, another bad region, or a region with distance $\leq d(D)-1$. See \cite{SW2008, OP2007} for details of the finger moves.
Then, we can reduce the number of edges in $\partial D$ by the above procedure. 

Repeat this operation inductively, and we get a nice diagram.
\end{proof}

Let $( \Sigma, \boldsymbol{ \beta }, \boldsymbol{ \alpha }, z_{0} )$ be a nice diagram of $-Y$ and $\phi \in \pi _{2} ^{0} ( \mathbf{x}, \mathbf{y} )$ be an empty embedded bigon. $-x$ and $y$ denote two vertices of $\phi$, with the signature of the orientation induced by $\partial ( \partial ( \phi ) | _{\boldsymbol { \beta } } )$. Then, there are $x_{2}, x_{3}, \ldots, x_{n}$ in intersection points between $\alpha$-curves and $\beta$-curves so that we can express $\mathbf{x} = x + x_{2} + x_{3} + \cdots + x_{n}$ and $\mathbf{y} = y + x_{2} + x_{3} + \cdots + x_{n}$. Such $x_{2}, x_{3}, \ldots, x_{n}$, $y$, and $x$ are called \textit{trivial vertices}, a \textit{positive vertex}, and a \textit{negative vertex} of $\phi$, respectively.

Similarly, when $\phi \in \pi _{2} ^{0} ( \mathbf{x}, \mathbf{y} )$ is an empty embedded square, and $-x, -x', y$, and $y'$ denote four vertices of $\phi$, there exist $x_{3}, x_{4}, \ldots, x_{n}$ in intersection points between $\alpha$-curves and $\beta$-curves so that we have an expression $\mathbf{x} = x + x' + x_{3} + x_{4} + \cdots + x_{n}$ and  $\mathbf{y} = y + y' + x_{3} + x_{4} + \cdots + x_{n}$. Such $x_{3}, x_{4}, \ldots, x_{n}$,  $y, y'$, and $x, x'$ are also called trivial vertices, positive vertices, and negative vertices of $\phi$, respectively.

The next lemma is useful for computing the boundary map $\hat{\partial}$.
\begin{Lem} \label{c-lemma}
Let $( \Sigma, \boldsymbol{ \beta }, \boldsymbol{ \alpha }, z_{0} )$ be a nice diagram of $-Y$ given by an open book.

Then, $c_{i}$ is not a negative vertex of any empty embedded bigon or square. Therefore, when there exists an empty embedded bigon or square $\phi$ in $\pi _{2} ^{0} ( \mathbf{x}, \mathbf{y} )$, if a generator $\mathbf{x}$ has a term $c_{i}$, $\mathbf{y}$ also does.

In particular, for any generator $\mathbf{x}$, $\pi _{2} ^{0} ( \mathbf{c}, \mathbf{x} )$ has no empty embedded bigons nor squares, where $\mathbf{c} = \sum c_{i}$. Hence, $\hat{\partial} \mathbf{c} = 0$ holds.
\end{Lem}
\begin{proof}
Let $\phi$ be an empty embedded bigon or square which has $c_{i}$ as a non-trivial vertex. Then, $\phi$ contains a thin strip between $a_{i} \times \{ 1/2 \}$ and $b_{i} \times \{ 1/2 \}$.
We orient $a_{i} \cap \partial \phi$ as the direction away from $c_{i}$. By the definition of $b_{i}$, $\phi$ is in the left side of $a_{i}$ as in Figure  \ref{zu-lem}.
Therefore, $c_{i}$ has positive signature.
\end{proof}
\myfig{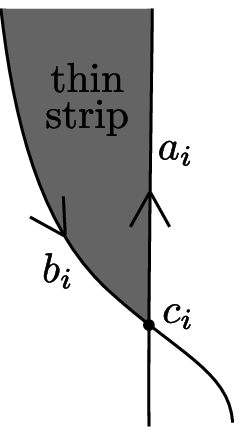}{Position of a thin strip.}{zu-lem}

\section{Examples of computing Ozsv\'ath-Szab\'o invariant} \label{exam}

In this section, we perform easier computation of the Ozsv\'ath-Szab\'o invariant for the open book which is investigated by Plamenevskaya in \cite{OP2007}.

Let $S$ be a four punctured sphere, and a monodromy $h$ be the composition of two maps $h = t_{ \gamma _{2} } \circ t _{ \gamma _{1} }$, where $t _{ \gamma _{1} }$ and $t_{ \gamma _{2} }$ are the positive Dehn twists around $\gamma _{1}$ and $\gamma _{2}$, respectively, as in Figure \ref{zu-exam-monod}.
\myfig{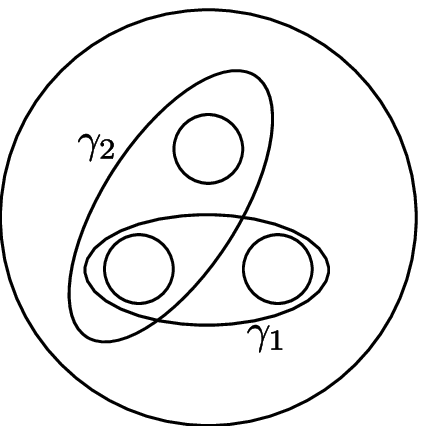}{Simple closed curves which give the definition of the monodromy $h$.}{zu-exam-monod}
It is known that the contact structure compatible with this open book is isotopic to the standard tight contact structure on $S^{1} \times S^{2}$. In fact, we can show, by combinatorial computation, that $c(\xi)$ is non-vanishing. 
\myfig{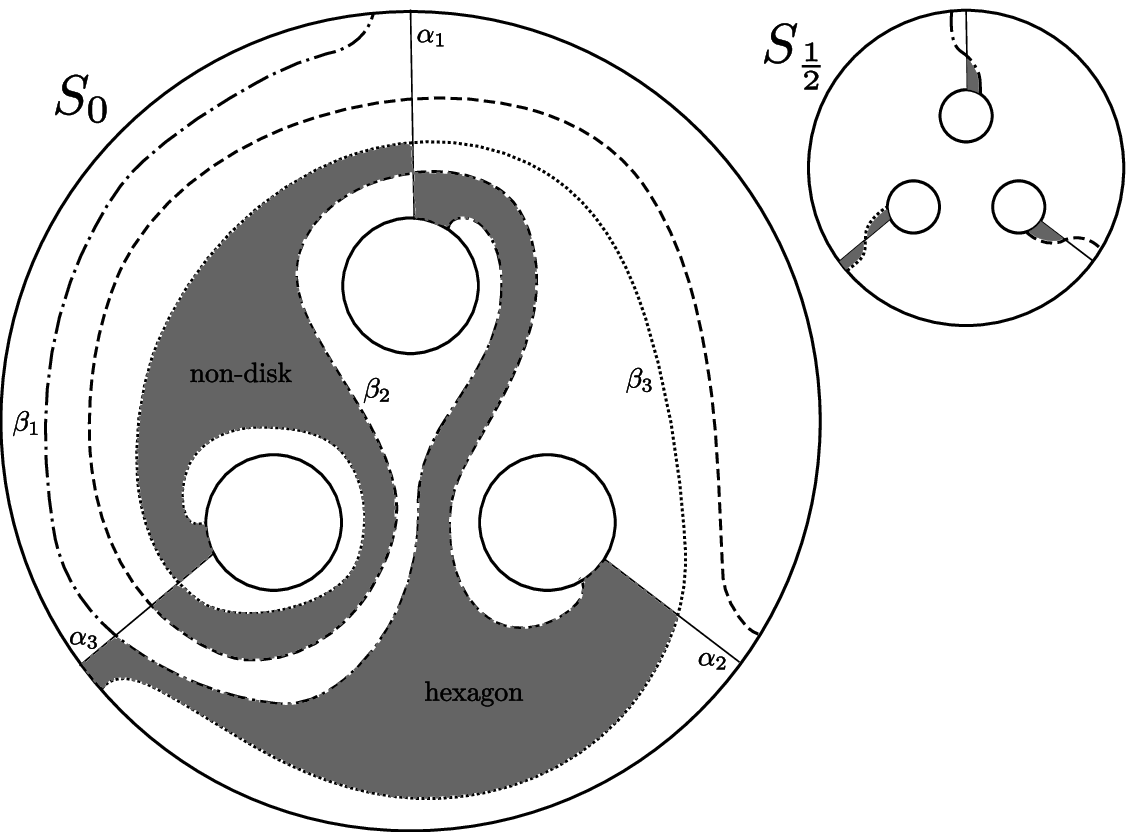}{Plamenevskaya's basis and two bad regions.}{zu-Pla-bad}

Plamenevskaya took arcs in Figure \ref{zu-Pla-bad} as a basis $\{ a_{1}, a_{2}, a_{3} \}$. The diagram given by her basis is shown by Figure \ref{zu-Pla-bad}, which has a non-disk region and a hexagon region.
\myfig{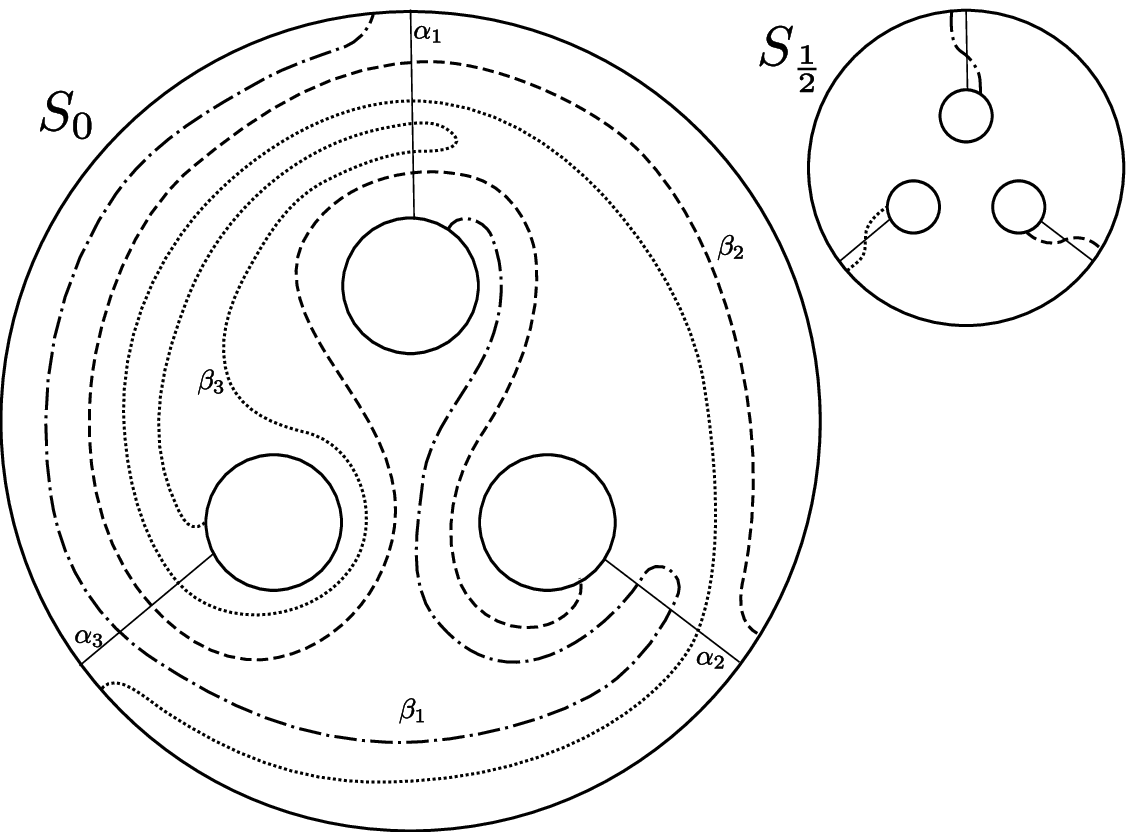}{The nice diagram after applying Sarkar-Wang algorithm to the diagram in Figure \ref{zu-Pla-bad}.}{zu-Pla-nice}
Figure \ref{zu-Pla-nice} shows the nice diagram after applying Sarkar-Wang algorithm, which has $22$ generators of $\widehat{CF}$ and $13$ empty embedded bigons/squares. See \cite{EO2007} for details and computation.

In \cite{EO2007}, Etg\"u and Ozbagci get a diagram having less generators and empty embedded bigons/squares by choosing a good basis $\{ a_{1}, a_{2}, a_{3} \}$. Here, we take another basis $\{ a_{1}, a_{2}, a_{3} \}$ as in Figure \ref{zu-exam-nice}.

In the corresponding diagram, $\alpha$-curves and $\beta$-curves have intersection points as shown in Table \ref{hyo-exam-points}.
 \begin{table}
  \caption{The intersection points between $\alpha$-curves and $\beta$-curves in Figure \ref{zu-exam-nice}. The components of $\alpha _{i}$-line and $\beta _{j}$-column mean intersection points between $\alpha _{i}$ and $\beta _{j}$.} 
  \label{hyo-exam-points}
   \begin{center}
    \begin{tabular}{|c|c|c|c|}
     \hline
                   & $\beta _{1}$ & $\beta _{2}$ & $\beta _{3}$ \\
     \hline
     $\alpha _{1}$ &    $c_{1}$   & $2$ points   & $2$ points   \\
     \hline
     $\alpha _{2}$ &  $\emptyset$ & $c_{2}$      & $x$          \\
     \hline
     $\alpha _{3}$ &  $\emptyset$ & $y$          & $c_{3}$      \\
     \hline
    \end{tabular}
   \end{center}
 \end{table}
Each generator of $\widehat{CF}$ corresponding to the above diagram has the term $c_{1}$ but does not contain any other terms included in $\alpha _{1}$, since $\beta _{1}$ intersects the $\alpha$-curves only in one point $c_{1}$, which is included in  $\alpha _{1}$. 
Thus, $\widehat{CF}$ has only $2$ generators, $\mathbf{c} = c_{1} + c_{2} + c_{3}$ and $\mathbf{x} = c_{1} + x + y$.

Moreover, there are only $2$ empty embedded squares for this diagram, $\phi _{1} = R_{1} + R_{2}$ and $\phi _{2} = R_{3} + R_{4}$, both of which are in $\pi _{2} ^{0} ( \mathbf{x}, \mathbf{c} )$, because the other square domains have a vertex which is an intersection point between $\alpha _{1}$ and $\beta _{2}$, or $\alpha _{1}$ and $\beta _{3}$.

So, we have $\hat{\partial} \mathbf{x} = \mathbf{c} + \mathbf{c} = 0$. Thus, 
\[ \widehat{ HF } ( S, h ) = \mathbb{Z} / 2 \mathbb{Z} [ \mathbf{c} ] \oplus \mathbb{Z} / 2 \mathbb{Z} [ \mathbf{x} ] \]
holds, and Ozsv\'ath-Szab\'o invariant $c(\xi) = [\mathbf{c}]$ does not vanish.
\myfig{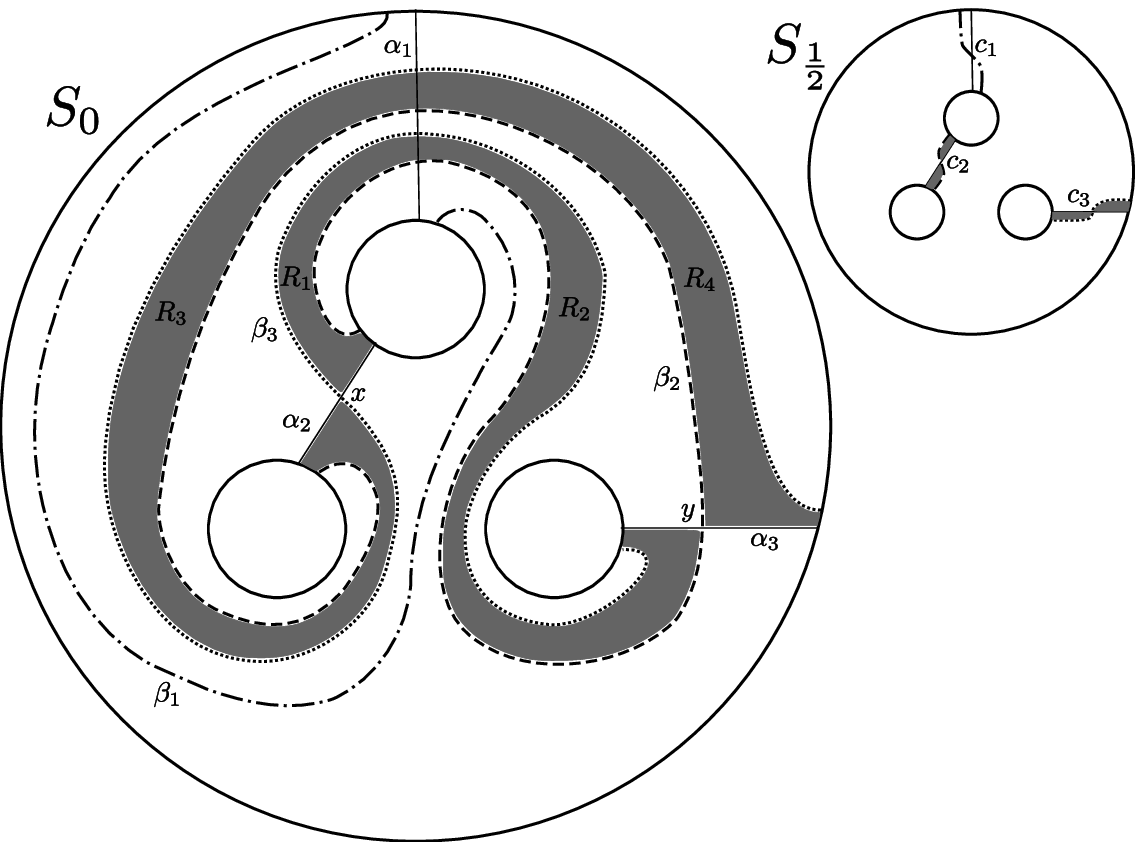}{One of good bases $\{ a_{1}, a_{2}, a_{3} \}$ and two empty embedded squares.}{zu-exam-nice}

\section{Ozsv\'ath-Szab\'o invariant of negative stabilizations} \label{nega-stab}

In this section, we show that the Ozsv\'ath-Szab\'o invariant for a contact structure $\xi$ vanishes, if $\xi$ is compatible with a negative stabilization of an open book.

Let $(S, h)$ be an open book and $(S', h')$ be its negative stabilization, i.e., $S'$ is a surface obtained by attaching a $1$-handle to $S$ and $\gamma$ is a simple closed curve in $S'$ which goes over the $1$-handle once and $h'$ is a composition map $h' = t_{\gamma} \circ h $, where $t_{\gamma}$ is the negative Dehn twist around $\gamma$. 

It is well-known that a contact structure which is compatible with $(S', h')$ is overtwisted. For a proof, see \cite{NG2003}, for example. So, the Ozsv\'ath-Szab\'o invariant vanishes. We show this combinatorially.

First of all, we take a basis $\{ a_{0}, a_{1}, \ldots, a_{n} \}$ for $S'$, such that $a_{0}$ is in the $1$-handle, which intersects $\gamma$ at only one point, and $\{ a_{1}, a_{2}, \ldots, a_{n} \}$ is a basis for $S$ which does not intersect $\gamma$.

In the diagram of $(S', h')$ given by the above basis, $S' _{1/2}$ includes an intersection point $c_{0}$ between $\beta _{0}$ and $\alpha _{0}$, and $S' _{0}$ includes two intersection points $x$, $y$ between $\beta _{0}$ and $\alpha _{0}$. Moreover, there are no intersection points between $\beta _{0}$ and the other $\alpha _{i}$'s.
\myfig{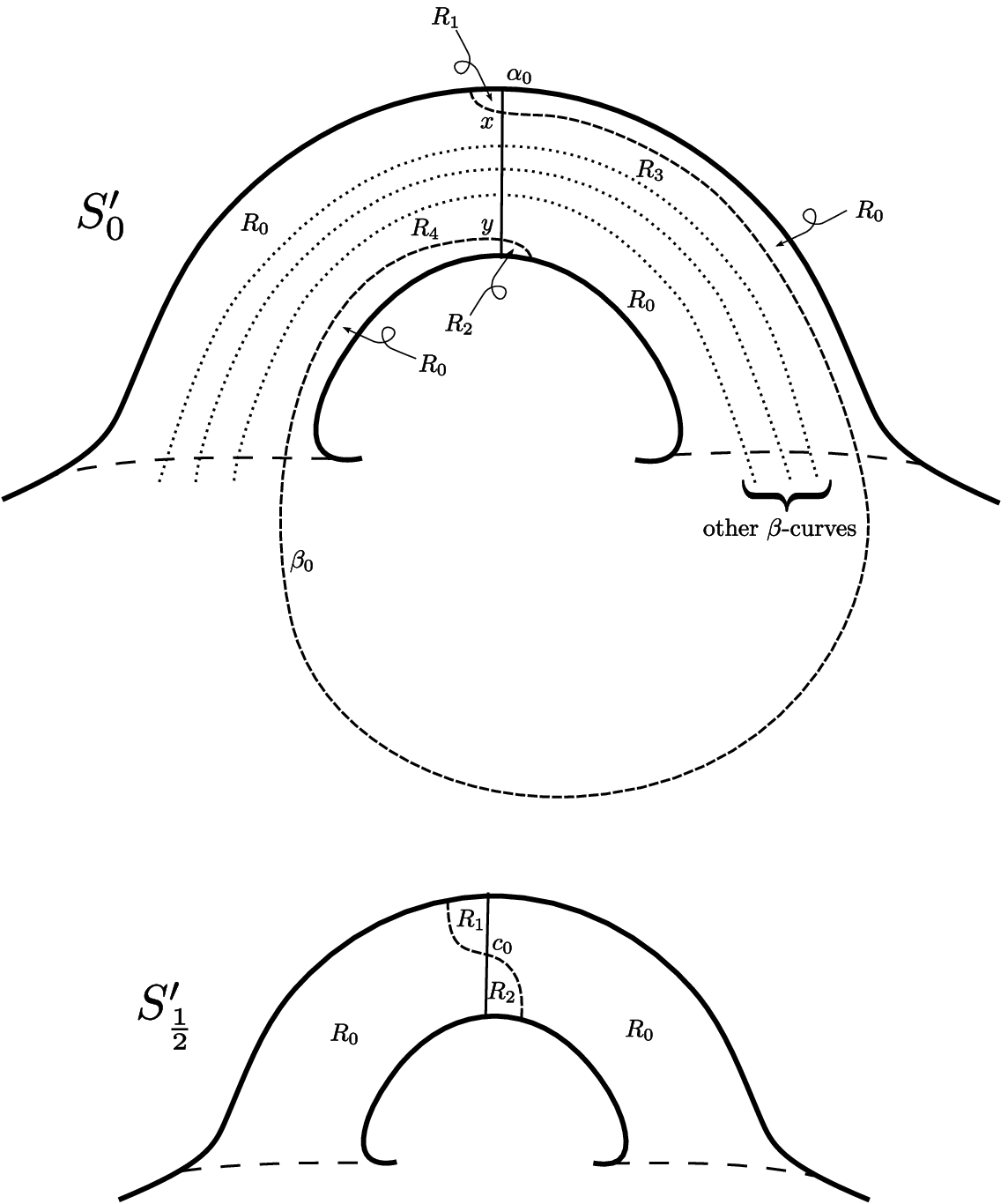}{The situation around $\beta _{0}$ and the $1$-handle of the negative stabilization.}{zu-nega-stab}

In Figure \ref{zu-nega-stab}, $R_{0}$ is the pointed region, and $R_{1}$, $R_{2}$, $R_{3}$, and $R_{4}$ are regions adjacent to $R_{0}$. But $R_{3}$ and $R_{4}$ can be connected to $R_{0}$ by arcs near $\beta _{0}$ which are disjoint from $\alpha$- and $\beta$- curves, because $\beta _{0}$ does not intersect $\alpha$-curves except $\alpha _{0}$. Hence, $R_{3}$ and $R_{4}$ coincide with $R_{0}$, and $\beta _{0}$ and bigon regions $R_{1}$, $R_{2}$ are surrounded with the pointed region. 

Therefore, $\beta _{0}$ and the bigon regions $R_{1}$, $R_{2}$ remain unchanged in a new nice diagram, which is obtained by applying Sarkar-Wang algorithm to the original diagram. 

Let $\mathbf{x} = x + c_{1} + c_{2} + \cdots + c_{n}$, $\mathbf{y} = y + c_{1} + c_{2} + \cdots + c_{n}$, and $\mathbf{c} = c_{0} + c_{1} + c_{2} + \cdots + c_{n}$ be the generators of $\widehat{CF}$ corresponding to the new nice diagram.
\begin{Thm}
We have
\[ \hat{\partial} \mathbf{x} = \hat{\partial} \mathbf{y} = \mathbf{c}. \]
Therefore, the Ozsv\'ath-Szab\'o invariant vanishes.
\end{Thm}
\begin{proof}
By lemma \ref{c-lemma}, if there exists a generator $\mathbf{x'}$ and an empty embedded bigon/square $\phi$ in $\pi _{2} ^{0} (\mathbf{x}, \mathbf{x'})$, then $\mathbf{x'}$ can be expressed as $x' + c_{1} + c_{2} + \cdots + c_{n}$, where $x'$ is an intersection point between $\alpha _{0}$ and $\beta _{0}$, and $\phi$ is a bigon with two vertices, $x'$ as a positive vertex and $x$ as a negative vertex.

Now, there are only three intersection points $c_{0}$, $x$, and $y$ between $\alpha _{0}$ and $\beta _{0}$. Clearly, there is no bigon whose vertices are $x$ and $y$. On the other hand, $R _{1}$ is the only empty embedded bigon with vertices $c_{0}$ and $x$.

Therefore, $\hat{\partial} \mathbf{x} = \mathbf{c}$

Similarly, we get $\hat{\partial} \mathbf{y} = \mathbf{c}$.
\end{proof}

The next theorem is Lemma 5.1.2 of \cite{NG2003}.
\begin{Thm}
If $\xi$ is an overtwisted contact structure on $Y$, then there exists an open book $(S, h)$ such that $\xi$ is compatible with a negative stabilization of $(S, h)$.
\end{Thm}

By using this theorem, we get the following.
\begin{Cor}
If $\xi$ is an overtwisted contact structure on $Y$, then $c(\xi)$ vanishes.
\end{Cor}

\section{Ozsv\'ath-Szab\'o invariant of positive stabilizations}
In this section, we construct an isomorphism between homology of an open book $(S, h)$ and that of the positive stabilization $(S', h')$ which preserves Ozsv\'ath-Szab\'o invariant.

Let $(S, h)$ be an open book and $(S', h')$ be its positive stabilization, i.e., $S'$ is a surface obtained by attaching a $1$-handle to $S$, and $\gamma$ is a simple closed curve in $S'$ which goes over the $1$-handle once, and $h'$ is a composition map $h' = t_{\gamma} \circ h $, where $t_{\gamma}$ is the positive Dehn twist around $\gamma$. 

For a contact $3$-manifold $(Y, \xi)$ and the compatible open book $(S, h)$, let $\xi ^{ \prime }$ denote the contact structure compatible with $(S', h')$, which is a positive stabilization of $(S, h)$. Then, it is a well-known fact that $\xi$ and $\xi ^{ \prime }$ are isotopic. For a proof, see \cite{NG2003}, for example.

We start the construction with taking a basis $\{ a_{0}, a_{1}, \ldots, a_{n} \}$ for $S'$ as in section \ref{nega-stab}, i.e., $a_{0}$ is in the $1$-handle, which intersects $\gamma$ at only one point, and $\{ a_{1}, a_{2}, \ldots, a_{n} \}$ is a basis for $S$ which is disjoint from $\gamma$. Then, by changing $h$, if we need, applying Theorem \ref{Pla-th}, we may assume that the diagram of $(S, h)$ given by the basis $\{ a_{1}, a_{2}, \ldots, a_{n} \}$ is nice.

Deforming $\gamma$ by isotopy in $S \setminus \left ( \bigcup_{i=1}^{n} a_{i} \right )$, we minimize the number of intersection points between $S \cap \gamma$ and $h(b_{1})$, $h(b_{2})$, \ldots, $h(b_{n})$. We may assume that the intersection of $\gamma \cap S$ and a bigon or a square region $R$ consists of arcs $\{ A_{i} \}$ in $R$, whose endpoints are on the $\beta$-edges contained in $h(b_{1}) \cup \cdots \cup h(b_{n})$. If at least one of such arcs $A_{i}$ has the endpoints on a common $\beta$-edge, we can choose $A_{i}$ so that $A_{i}$ and the $\beta$-edge bound a bigon whose interior disjoint from $\gamma \cap S$. Then deforming $\gamma$ by isotopy, we can eliminate $A_{i}$. This contradicts the above minimality condition. It follows that $\gamma \cap S$ does not intersect a bigon, and intersects a square region in arcs, each of which connects the different $\beta$-edges of the region. Note that this deforming does not change the equivalence class of the original open book.
\myfig{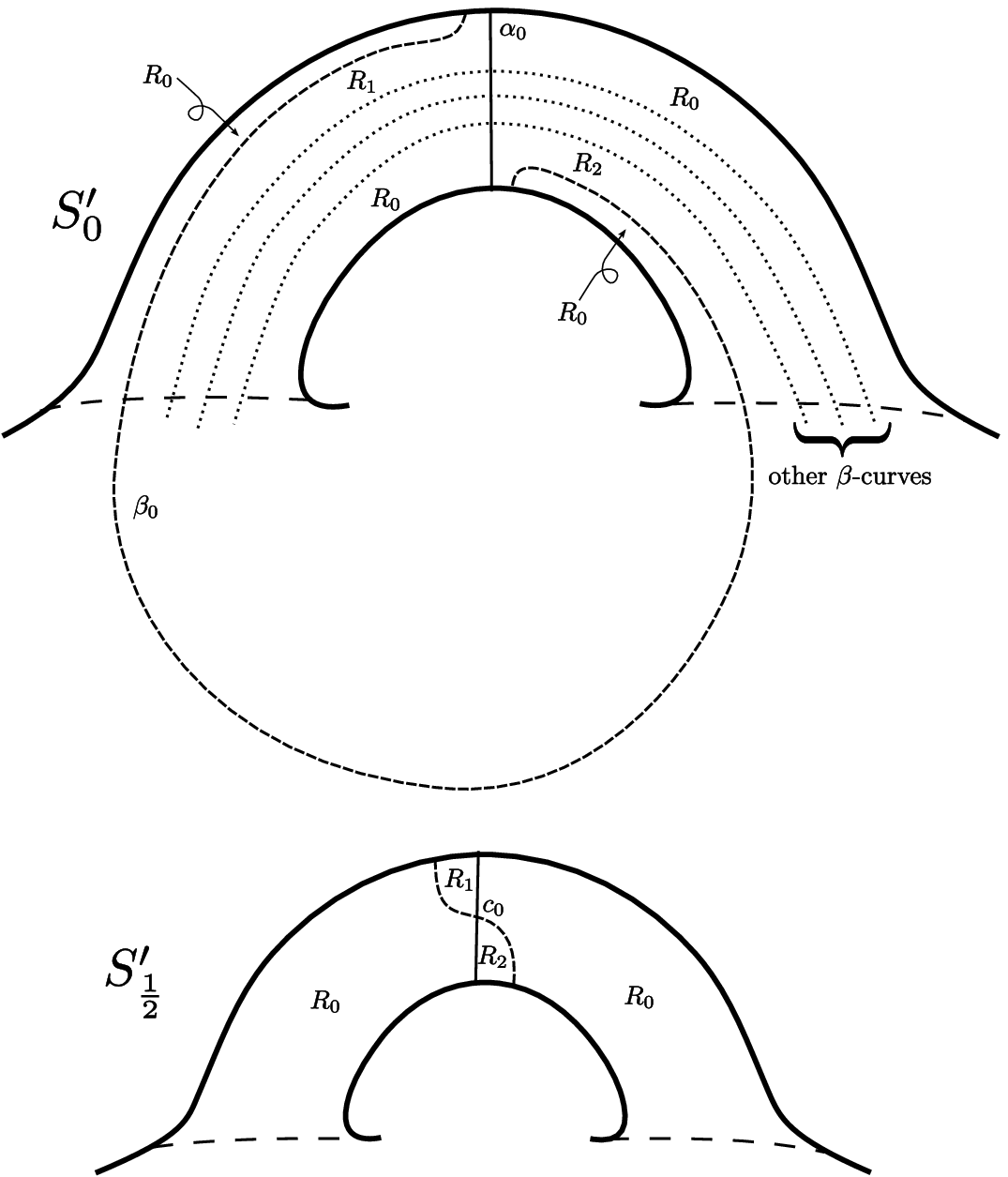}{The situation around $\beta _{0}$ and the $1$-handle of the positive stabilization.}{zu-posi-stab}

In the diagram of $(S', h')$ given by the above basis, $\beta _{0}$ intersects the $\alpha$-curves only at one point $c_{0} = \alpha _{0} \cap \beta _{0}$. So, the same argument as in the case of negative stabilization shows that only the pointed region has an edge contained in $\beta _{0}$. In fact, in Figure \ref{zu-posi-stab} a region $R_{1}$ can be connected to the pointed region $R_{0}$ by an arc near $\beta _{0}$ and a region $R_{2}$ can be connected to $R_{0}$ similarly. Therefore, $R_{1}$ and $R_{2}$ coincide with the pointed region $R_{0}$. Hence, non-pointed regions have their boundaries in $\alpha _{0}, \alpha _{1}, \ldots, \alpha _{n}$ and $\beta _{1}, \beta_{2}, \ldots, \beta _{n}$.

We prove that this diagram is nice. To see this, we take a non-pointed region $R$. If the boundary of $R$ does not intersect $\alpha _{0}$, then $R$ coincides with a region of diagram for $(S, h)$, which remains unchanged by the stabilization, so $R$ is a bigon or a square.

On the other hand, if the boundary of $R$ does intersect $\alpha _{0}$, then $R$ is obtained from some region $\widetilde{R}$ of the diagram for $(S, h)$ by Dehn twist around $\gamma$ as in Figure \ref{zu-twist}. Hence, $\widetilde{R}$ is a square region which intersects $\gamma$. The stabilization changes the regions only around a tubular neighborhood of $\gamma$ as in Figure \ref{zu-tubularnbd}. By the stabilization, each square region which intersects $\gamma$ in the original diagram is twisted and is divided into several regions by $\alpha _{0}$. Therefore, the resulting new regions are still square regions. 
\myfig{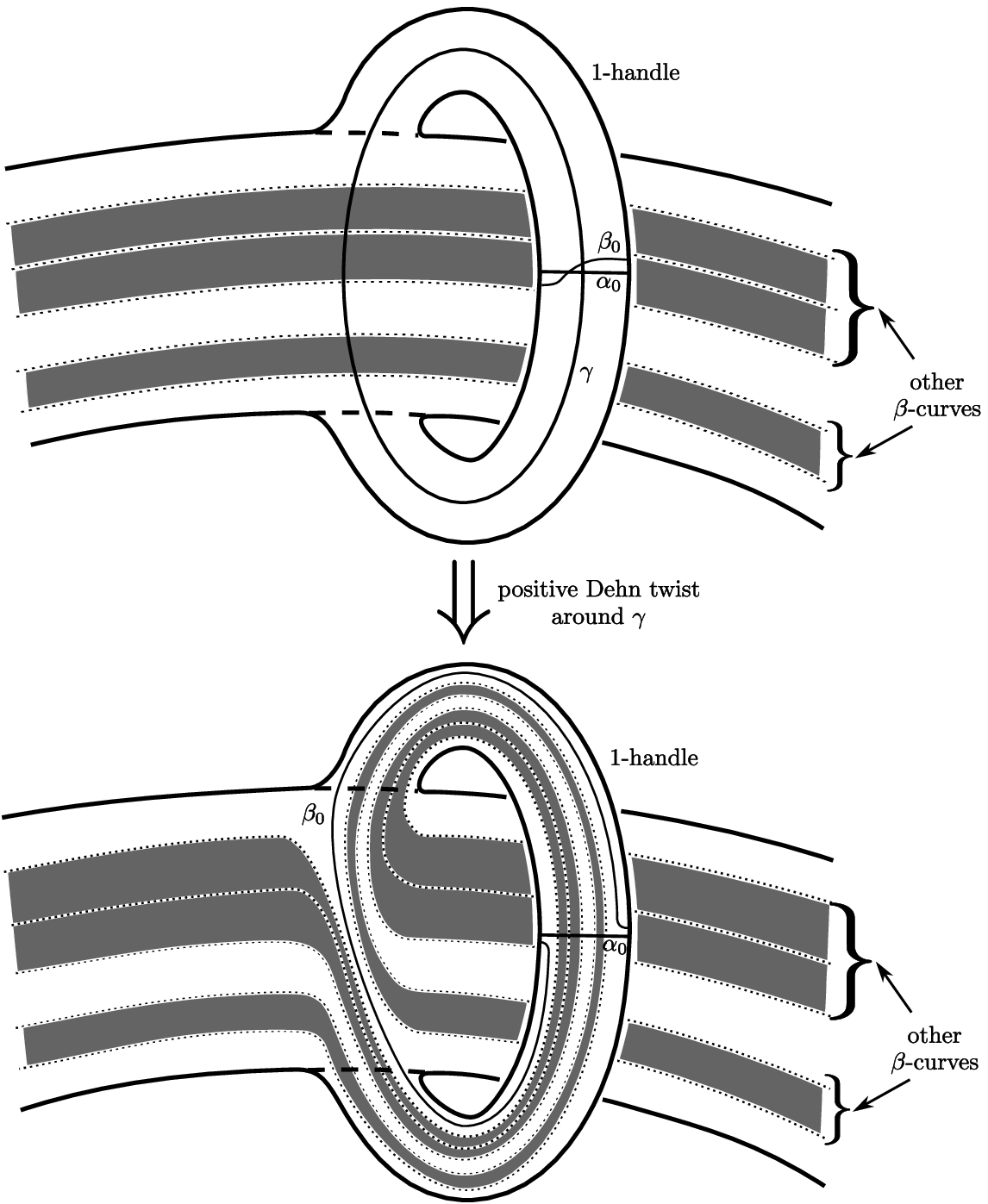}{The correspondence between empty embedded bigons or squares.}{zu-twist}
\myfig{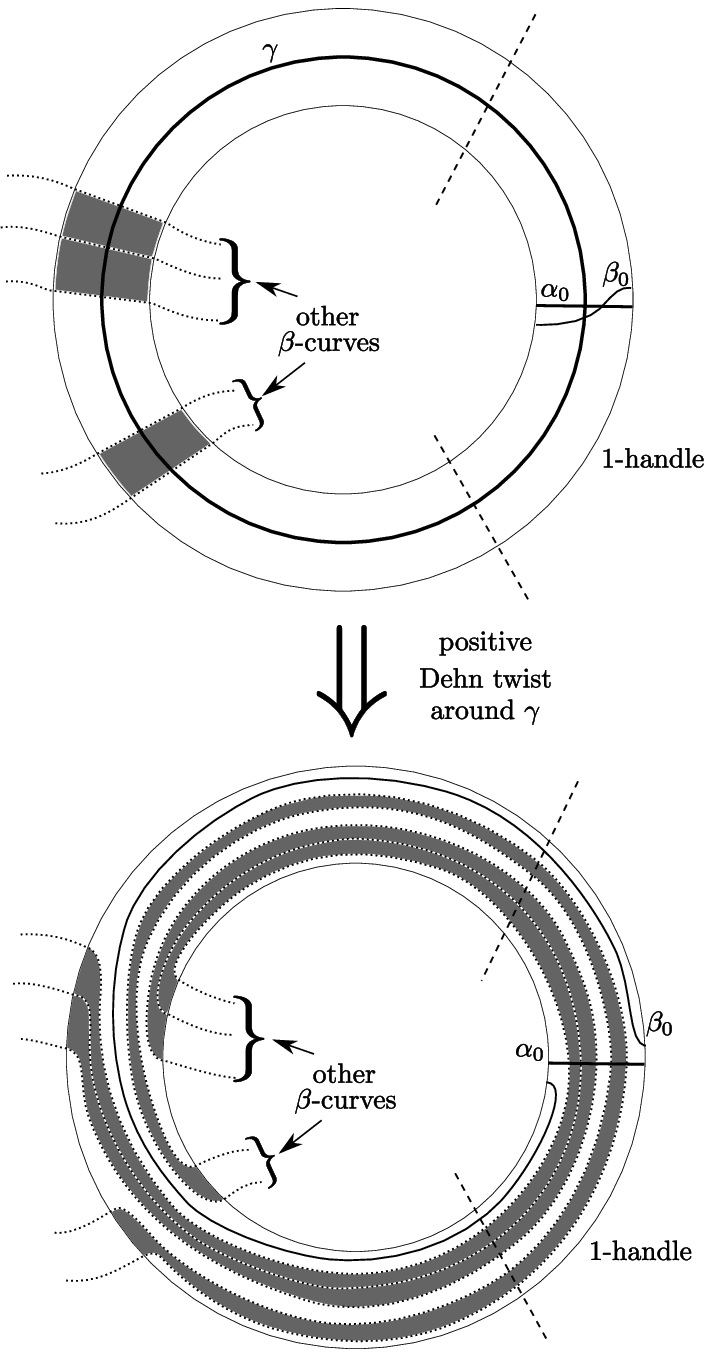}{How positive stabilization changes curves and regions in the tubular neighborhood around $\gamma$.}{zu-tubularnbd}

Each generator $\bar{ \mathbf{x} }$ of $\widehat{CF}(S',h')$ has the term $c_{0}$, since $\beta _{0}$ intersects the $\alpha$-curves only at $c_{0}$, which is included in $\alpha _{0}$. The other terms of $\bar{ \mathbf{x} }$ are intersection points between $\alpha _{1}, \alpha _{2}, \ldots, \alpha _{n}$ and $\beta _{1}, \beta_{2}, \ldots, \beta _{n}$. Therefore, we get a one-to-one correspondence between $\widehat{ CF } (S, h)$ and $\widehat{ CF } (S', h')$ by mapping a generator $\mathbf{x}$ of $\widehat{ CF } (S, h)$ to $\bar{ \mathbf{x} } = c_{0} + \mathbf{x}$.

Additionally, for an empty embedded bigon or square $\phi \in \pi _{2}^{0} (\mathbf{x}, \mathbf{y})$ in the diagram for $(S, h)$, we assign to $\phi$ an empty embedded bigon or square $\bar{ \phi } \in \pi _{2}^{0} (\bar{ \mathbf{x} }, \bar{ \mathbf{y} })$ which goes through the $1$-handle as in Figure \ref{zu-twist}, if $\phi$ intersects $\gamma$ in the diagram for $(S, h)$, and assign $\phi$ itself to $\phi$, otherwise(note that $\phi$ remains unchanged by stabilization in the latter case).
Then this correspondence is one-to-one between the empty embedded bigons or squares connecting $\mathbf{x}$ to $\mathbf{y}$ and the empty embedded bigons or squares connecting  $\bar{ \mathbf{x} }$ to $\bar{ \mathbf{y} }$.

After all, we get the following theorem.

\begin{Thm}
There exists an isomorphism $i$ between $\widehat{ CF } (S, h)$ and $\widehat{ CF } (S', h')$ such that each generator $\mathbf{x}$ is mapped to $\bar{ \mathbf{x} } = c_{0} + \mathbf{x}$.

In addition, there is a one-to-one correspondence between the empty embedded bigons or squares connecting $\mathbf{x}$ to $\mathbf{y}$ and the empty embedded bigons or squares connecting  $\bar{ \mathbf{x} }$ to $\bar{ \mathbf{y} }$, for all generators $\mathbf{x}, \mathbf{y}$ of $\widehat{ CF } (S, h)$.

Therefore, $i$ is an isomorphism between the chain complexes, and the induced isomorphism $i_{\ast} \colon \widehat{ HF } (S, h) \rightarrow \widehat{ HF } (S', h')$ satisfies $i_{\ast} ( c( \xi ) ) = c( \xi ^{\prime} )$.
\end{Thm}

\bibliographystyle{plain}
\bibliography{ref-arXiv}

\end{document}